\documentclass[11pt]{amsart}
\usepackage{amscd,bbm,dsfont}
\usepackage{amsmath,amssymb,mathrsfs,amsthm}
\usepackage{amsfonts}
\theoremstyle{plain}
\newtheorem{theorem}{Theorem}[section]
\newtheorem{corollary}[theorem]{Corollary}
\newtheorem{lemma}[theorem]{Lemma}
\newtheorem{proposition}[theorem]{Proposition}
\theoremstyle{definition}
\newtheorem{definition}{Definition}[section]
\theoremstyle{remark}
\newtheorem{remark}{Remark}[section]
\numberwithin{equation}{section}

\begin{document}
\title[A class of degenerate elliptic equations]
{\Small Well-posedness of boundary value problems for a class of
degenerate elliptic equations}
%\vskip5pt
\author{Yue HE}
%Address of record for the research reported here
\address{address 1: Department of Mathematics, School of Mathematics \& Computing,
Nanjing Normal University, Nanjing 210097, P. R. China}
%Current address
\curraddr{Department of Mathematics, Zhongshan University, Guangzhou
510275, P. R. China} \email{heyueyn@163.com \& heyue@njnu.edu.cn.}
\thanks{Partially supported by NNSF grant of P. R. China: No.10571087,
and Natural Science Foundation of Jiangsu Education Commission of P.
R. China: No.05KJB110063}
% ------------------------------------------------------------------------
% ------------------------------------------------------------------------
%\thanks{This work was completed with the
%support of NNSF of China: No.10571087, and Natural Science
%Foundation of Jiangsu Education Commission: No.05KJB110063}
%\thanks{Project(No.04KJB110062) supported in part by the
%Natural Science Foundation of Jiangsu Province of P. R. China}
% ------------------------------------------------------------------------
% ------------------------------------------------------------------------
%General information
\subjclass[2000]{Primary 35J70; Secondary 35B65, 35J60, 35J25.}
\keywords{rigidity problem arising in infinitesimal isometric
deformation, degenerate elliptic equations, prior estimate,
well-posedness.}
%\author{\small School of Mathematics \& Computing}
% ------------------------------------------------------------------------
\date{}
\maketitle
% ------------------------------------------------------------------------
%\baselineskip=14pt
% ------------------------------------------------------------------------
%\thanks{This work was completed with the support of .}
%\thanks{The author was also supported in part by CNNSF.}
%\thanks{This work was also supported in part by JSPNSF: 04KJB110062.}
%\date{July 28, 2004 and, in revised form, July 6, 2005.}

\begin{abstract}
In this paper, we study the well-posedness of boundary value
problems for a special class of degenerate elliptic equations coming
from geometry. Such problems is intimately tied to rigidity problem
arising in infinitesimal isometric deformation, The characteristic
form of this class of equations is changing its signs in the domain.
Therefore the well-posedness of these above problems deserve to make
a further discussion. Finally, we get the existence and uniqueness
of $H^1$ solution for such boundary value problems.
\end{abstract}
% ------------------------------------------------------------------------
% ------------------------------------------------------------------------
%\\{\bf Chinese Library Classification}\quad O175.2\quad\quad
% ------------------------------------------------------------------------
%\\{\bf Document Code A}\quad\quad
% ------------------------------------------------------------------------
%\\{\bf Article ID}\quad

\section{Introduction.}
\setcounter{equation}{0} In this section, we introduce briefly the
history of the research to degenerate elliptic equations and some
geometric backgrounds to a class of degenerate elliptic equations
which we are concerned with. Next we will put forward the main
question of this paper. In the last part of this section, we will
summarize the main result and its trivial generalization about
existence and uniqueness for solution to such a class of equations.

\subsection{Historical remarks and backgrounds.}
In this subsection, we will introduce the brief history and the
current status of investigation about the degenerate elliptic
equations. We observe the equation
$$Lu\equiv a^{ij}u_{x_ix_j}+
b^{k}u_{x_k}+cu=f\quad\hbox{in}\,\,\,\Omega,$$ where
$u_{x_k}=\partial u/\partial x_k, u_{x_ix_j}=\partial^2u/\partial
x_i\partial x_j$ etc. and the index $i,j,k$ runs from $1$ to $n$ and
repeated indices imply summation. \,\,\,(From now on we will use
such summation convention throughout this paper).\,\,\, If for any
vector $\xi=(\xi_{1},...,\xi_{n})\in \mathbb{R}^n$,
$a^{ij}(x)\xi_{i}\xi_{j}\geq 0$\,\, for all\,\,
$x=(x_{1},...,x_{n})\in \Omega\subseteq\mathbb{R}^n$. Then $Lu=f$ is
called second order PDE with nonnegative characteristic form, and
also called second order degenerate elliptic equation, or second
order elliptic-parabolic equation in domain $\Omega$. They contain
elliptic equation, parabolic equation, one order differential
equation, Brown Motion equation and some important equations
introduced later.

The study of degenerate elliptic equations can be traced back to
1910's, first appeared in Picone's thesis. After that, the Tricomi's
research report \cite{Tr} and his research on the mixed-type partial
differential equation. Besides, M.V. Keldy$\check{s}$ \cite{Ke};
Fichera; J. Kohn and L. Nirenberg \cite{KN}; O.A.Ole$\check{i}$nik
\cite{Ol}, \cite{Ole}; E.V. Radkevi$\check{c}$ \cite{Ra} etc., they
all do many works for laying a foundation in this field. After half
century's development in this field, O.A. Ole$\check{i}$nik and E.V.
Radkevi$\check{c}$, published their classical monograph \cite{OR} in
1971. Their monograph summarized the theories developed before
1970's, and established a general framework for the theories to
second order degenerate linear elliptic equations. They stated the
existence and uniqueness of weak solution for the general boundary
value problem to such equations in $L^p$ space and some Hilbert
spaces, and made a certain further contribution to the regularity
theories of weak solutions for the general boundary value problems
to such equations. During the past three decades, several progresses
have been made in the research of second order degenerate linear
elliptic equations. But we only enumerate the partial well-known and
representative work at here. For instance, L. Caffarelli, L.
Nirenberg and J. Spruck \cite{CNS} had studied the degenerate
Monge--Amp$\acute{e}$re equation; H. Brezis and P.L. Lions \cite{BL}
had studied the Yang--Mills equation $-x^{2}\Delta u+2u=f(u)$
describing gauge fields; E.B. Fabes, C.E. Kenig and D. Jerison
{\cite{FKJ} had studied the general degenerate elliptic equation
with the divergence form $\partial_{i}(a_{ij}(x)\partial_{j}u)=f(x)$
and so on.

\textbf{A geometric background.} Fanghua Lin \cite{Lin} had studied
the following Dirichlet problem for minimal graphs in hyperbolic
space
\begin{equation}\label{isothermal}
\left\{
\begin{array}{lll}
\Delta f-\displaystyle\frac{f_if_j}{1+|df|^2}f_{ij}+\frac{n}{f}=0
\quad \hbox{in}\,\,\,\Omega,\\[5pt]
f>0\quad\,\hbox{in}\,\,\,\Omega, \\
f\big|_{\partial\Omega}=0\quad\,\hbox{on}\,\,\,\partial\Omega,
\end{array}
\right.
\end{equation}
where $\Omega\subset\mathbb{R}^n,\,n\geq 2,$ is a bounded open
domain; $|df|^2=\Sigma_{i=1}^n\,f_i^2,\,\,f_i=f_{x_i}$ for $0\leq
i\leq n$. Obviously, by direct computation we know that the equation
\eqref{isothermal} is just the Euler-Lagrange equation of the
variational integral
$$A[f,K]=\int_K f^{-n}\sqrt{1+|df|^2}dx$$
for all compact subsets $K$ of $\Omega$.
Since $graph(f)$ is a $C^{1,\alpha}$-manifold with boundary in $\mathbb{R}^{n+1}$,
and since the tangent planes of $graph(f)$ along the boundary $\partial\Omega$ are vertical,
we view $graph(f)$ near a point at $\partial\Omega$ as a graph over
such a vertical plane.
This is equivalent to the hodograph transformation of \eqref{isothermal}.
Then the corresponding P.D.E. for the function $u$ which represents the $graph(f)$ is
\begin{equation}\label{heinz}
\left\{
\begin{array}{ll}
y(\Delta
u-\displaystyle\frac{u_{i}u_{j}}{1+|du|^2}u_{ij})-nu_{y}=0\quad
\hbox{in}\,\,\,B_{1}^{+}(0),\\[5pt]
u(x,0)=\varphi(x)\quad \hbox{for}\,\,\,|x|\leq 1,
\end{array}
\right.
\end{equation}
where $B^{+}_{1}(0)=\{(x,y)\in \mathbb{R}^{n}_{+}:|x|\leq 1,0\leq
y\leq 1 \}$. This is a degenerate quasilinear elliptic equation.
Fanghua Lin shows the solutions $u$ of \eqref{heinz} are as smooth
as $\varphi$ in $B^{+}_{1/2}(0)$, and use this fact to prove the
following result:

\textbf{Theorem.}\,\,\emph{If $\partial\Omega$ is of class
$C^{k,\alpha}$, then $graph(f)$ is a $C^{k,\alpha}$ hypersurface
with boundary for either (1)\,\,$1\leq k\leq n-1$ and $0\leq
\alpha\leq 1$ or (2)\,\,$n\leq k\leq \infty$ and $0<\alpha<1$.}

\textbf{Another geometric background.}
Recently, ones come into
contact with a class of second order degenerate elliptic equations
when they study the rigidity problem arising in infinitesimal
isometric deformation. We shall simply introduce some geometric
backgrounds about the above equations in the following. The details
can be found in \cite{STY}, \cite{Hong}, \cite{HJX} and
\cite{Lihong}.

Given a metric $g$ with smooth positive curvature $K$ on the
closed unit disk $\bar{D}$. In the sequel we always denote it by
$(\bar{D},g)$. In terms of local coordinates system $(u^1,u^2)$ on
$\bar{D}$, the metric $g$ can be expressed as $g=g_{ij}du^idu^j$.
Suppose that $\vec{r}=(x,y,z)$ is a smooth isometric embedding of
$(\bar{D},g)$ into $\mathbb{R}^3$. and the boundary
$\vec{r}(\partial D)$ is a $C^2$ planar  convex curve. By the
Gauss equations we have in a local coordinate system,
\begin{equation}\label{sobolov}
\vec{r}_{ij}=\Gamma_{ij}^k \vec{r}_k+\Omega_{ij}\vec{n}
\,\,\,\hbox{or}\,\,\,\nabla_{ij}\vec{r}=\Omega_{ij} \vec{n}\,\,\quad(i,j,k=1,2),
\end{equation}
where subscripts $i,j$ and $\nabla_{ij}$ denote Euclidean and
convariant derivatives respectively, $\Omega_{ij}$ the
coefficients of the second fundamental form, $\Gamma_{ij}^k$ the
Christoffel symbols with respect to the metric and $\vec{n}$ the
unit normal to $\vec{r}$. For each unit constant vector, for
instance, the unit vector $\vec{k}$ of the $z$ axis, taking the
scale product of $\vec{k}$ with the two hind sides of
\eqref{sobolov} and using the Gauss equations one can get
\begin{equation}\label{janet}
\det(\nabla_{ij}z)=K\det(g_{ij})\,(\vec{n},\vec{k})^2\quad\,(i,j=1,2),
\end{equation}
where $K$ is Gaussian curvature. Notice that
$$(\vec{n},\vec{k})^2
=1-\Big(\frac{(\vec{r}_1\times\vec{r}_2)\times\vec{k}}{|\vec{r}_1\times\vec{r}_2|}\Big)^2
=1-g^{ij}z_iz_j =1-|\nabla z|^2,$$ where $\nabla
z=(g^{1l}z_l,g^{2l}z_l)$ is the gradient of $z$. Inserting the
last expression into \eqref{janet}, we deduce the Darboux equation
\begin{eqnarray}\label{burstin}
F(z)=\det(\nabla_{ij}z)-K\,\det(g_{ij})\,(1-|\nabla z|^2)=0.
\end{eqnarray}
Obviously each component of $\vec{r}$ satisfies the Darboux equation \eqref{burstin}.

Given a smooth surface $\vec{r}$ in $\mathbb{R}^3$ one consider
its deformation $\vec{r}_t:(-\varepsilon,\varepsilon)\ni
t\rightarrow \mathbb{R}^3$ with $\vec{r}_0=\vec{r}$. If $t=0$ is a
critical point of the metric $g(t)=d\vec{r}^2_t$, we say that the
derivative with respect to $t$ of $\vec{r}_t$ at $t=0$ give rise a
first order infinitesimal isometric deformation of $\vec{r}$.
Denoting this infinitesimal isometric deformation by
$\vec{\tau}=(d\vec{r}_t/dt)(0)$, and we call it the first order
infinitesimal deformation vector, or the first order deformation
vector. So we have
\begin{equation}\label{stokers}
\frac{d}{dt}(d\vec{r}_t^2)\big|_{t=0}=(d\vec{r},d\vec{\tau})=0.
\end{equation}
Obviously, any rigid body motion of $\vec{r}$,
$\vec{\tau}=\vec{A}\times\vec{r}+\vec{B}$ for arbitrary two constant vectors $\vec{A}$ and $\vec{B}$,
is always a solution of \eqref{stokers} and such solutions are called trivial ones.
We say that $\vec{r}$ is of
infinitesimal rigidity for first order isometric deformation
if \eqref{stokers} has no nontrivial solution.
As is well known, for closed surface now we only know that
closed $C^2$ convex surfaces are infinitesimally rigid.
For a surface $\vec{r}$ with boundary, usually it is not infinitesimally rigid if
there is no restriction to the deformation on the the boundary of $\vec{r}$.
Therefore we must impose some condition, for instance,
$$(\vec{\tau},\vec{k})=0\,\,\,\hbox{on}\,\,\,\partial D,$$
where $\vec{k}$ is the unit vector of $z$ axis.

Let us consider an infinitesimal isometric deformation of surface $\vec{r}$,
$\vec{r}_\epsilon=\vec{r}+\epsilon\vec{\tau}$ where $\vec{\tau}=(\xi,\eta,\zeta)$
satisfies \eqref{stokers}.
Notice that $\epsilon=0$ is the critical point of
$g_\epsilon=d\vec{r}_\epsilon^2$ and hence, the diferentiation of
its Gaussian curvature $K(\epsilon)$
and Christoffel symbols $\Gamma_{ij}^k(\epsilon)\,\,\,(i,j,k=1,2)$
(i.e. connection coefficients) in $\epsilon$
are equal to zero at $\epsilon=0$.
Then differntiation of the Darboux equation \eqref{burstin} for $z+\epsilon\zeta$
with respect to $\epsilon$, letting $\epsilon=0$, gives
\begin{equation}\label{ambient}
F^{ij}(z)\,\nabla_{ij}\zeta+2K\,\det(g_{ij})\,(\nabla z,\nabla
\zeta)=0,
\end{equation}
where $F^{ij}(z)=\partial\det(\nabla_{kl})/\partial\nabla_{ij}z\,\,\,(i,j=1,2)$
is the algebraic cofactors of $\nabla_{ij}z$.

Since
$\nabla_{ij}\vec{r}=\vec{r}_{ij}-\Gamma_{ij}^k\vec{r}_k=\Omega_{ij}\vec{n}\,\,\,(i,j=1,2)$,
we have
$\nabla_{ij}z=z_{ij}-\Gamma_{ij}^kz_k=\Omega_{ij}(\vec{n},\vec{k})\,\,\,(i,j=1,2)$.
So \eqref{ambient} can be written as follows
\begin{equation}\label{simply}
(\vec{n},\vec{k})\Omega^{ij}\,\nabla_{ij}\zeta+2(\nabla z,\nabla \zeta)=0.
\end{equation}
where $(\Omega^{ij})=(\Omega_{ij})^{-1}$ is the inverse of matrix $(\Omega_{ij})$.

Obviously, the $(\Omega^{ij})$ is a positive definite matrix if
Gaussian curvature $K$ is positive. \eqref{simply} is of nonnegative
characteristic form in the subdomain $\{\,(u^1,u^2)\in
\bar{D}:(\vec{n},\vec{k})\geq 0\}$, is of nonpositive characteristic
form in the subdomain $\{\,(u^1,u^2)\in
\bar{D}:(\vec{n},\vec{k})\leq 0\}$, and is one order PDE in the
subdomain $\{\,(u^1,u^2)\in \bar{D}:(\vec{n},\vec{k})=0\}$.
Therefore, \eqref{simply} is characteristic degenerate, and its
characteristic form changing sign in domain $\bar{D}$.

The spherical crown is an example of this aspect (see \cite{Hong}).
Let us consider a spherical crown $\sum_\lambda=\{x^2+y^2+z^2=1:z\leq 0\}$.
In spherical coordinates
$$\Sigma_\lambda=\{\,(\sin\theta \cos\phi,\sin\theta\sin\phi,-\cos\theta)
\,\,\big|\,\,0\leq \phi\leq 2\pi,0\leq \theta \leq
\theta_*=\arccos(-\lambda)\,\},$$ where $\lambda$ is a positive
constant, and $\theta=0$ stand for the South pole.
$\Sigma_\lambda$ is the isometric embedding of the metric
$g=d\theta^2+sin^2\theta d\phi^2, 0\leq \theta \leq \theta_*$.
Since $\lambda>0$, so $\Sigma_\lambda$ contains the below
hemisphere. In the present case, \eqref{simply} may be written as
follows
\begin{equation}\label{wholee}
\cos\theta\big[(\sin\theta\zeta_\theta)_\theta
+\big(\frac{\zeta_\phi}{\sin\theta}\big)_\phi\big]
+2\sin^2\theta\zeta_\theta=0,\,\,\,\theta\in (0,\theta_*).
\end{equation}
with the constraint condition in here as follows:
$$\zeta=0\,\,\,\hbox{on}\,\,\,\theta=\theta_*\,\,\,
\hbox{and}\,\,\,\zeta\,\,\,\hbox{is bounded near}\,\,\,\theta=0.$$
Evidently \eqref{wholee} is elliptic as $\theta\neq\pi/2$. We
ought to show that: if $\theta_*>0$, then \eqref{wholee} is of
nonnegative characteristic form as $0\leq\theta\leq\frac{\pi}{2}$,
and is of nonpositive characteristic form as
$\frac{\pi}{2}\leq\theta\leq\theta_*$.

In the same way, one may consider
another infinitesimal isometric deformation of surface $\vec{r}$ in the following
$\vec{r}_\epsilon=\vec{r}+\epsilon\vec{\tau}_1+\epsilon^2\vec{\tau}_2+\cdots$.
Denote $\vec{\tau_1}=(\xi_1,\eta_1,\zeta_1)$, $\vec{\tau}_2=(\xi_1,\eta_2,\zeta_2)$, and so on.
We respectively call $\vec{\tau_1}, \vec{\tau}_2$ the first order deformation vector,
the second order deformation vector, etc.
Obviously,
$$g_\epsilon=(d\vec{r}_\epsilon,d\vec{r}_\epsilon)
=(d\vec{r},d\vec{r})+2\epsilon(d\vec{r},d\vec{\tau}_1)
+\epsilon^2[2(d\vec{r},d\vec{\tau}_2)+(d\tau_1,d\tau_1)]+O(\epsilon^3).$$
If $g_\epsilon=(d\vec{r},d\vec{r})+O(\epsilon^3)$, Then
$\vec{\tau_1}$ and $\vec{\tau}_2$ should satisfy the following
systems
$$\left\{
\begin{array}{ll}
(d\vec{r},d\vec{\tau}_1)=0,\\
(d\vec{r},d\vec{\tau}_2)=-\frac{1}{2}(d\tau_1,d\tau_1).
\end{array}
\right.$$ We may analogously definite and discuss the rigidity
 of
the second order, even higher order infinitesimal isometric
deformation of surface $\vec{r}$. But in here we only show that
$g_\epsilon$ gives rise to a second order infinitesimal isometric
deformation of $g$ if $g_\epsilon$ is equal to $g$ up to second
order, i.e.
$$(dg_\epsilon/d\epsilon)\big|_{\epsilon=0}
=(d^2g_\epsilon/d\epsilon^2)\big|_{\epsilon=0}=0.$$
Consequently,
$$\frac{d\,\Gamma_{ij}^k(\epsilon)}{d\epsilon}\big|_{\epsilon=0}
=\frac{d^2\Gamma_{ij}^k(\epsilon)}{d\epsilon^2}\big|_{\epsilon=0}
=\frac{dK(\epsilon)}{d\epsilon}\big|_{\epsilon=0}
=\frac{d^2K(\epsilon)}{d\epsilon^2}\big|_{\epsilon=0}=0.$$ And
then, Then two order derivative of the Darboux equation
\eqref{burstin} for $z+\epsilon\zeta$ with respect to $\epsilon$,
letting $\epsilon=0$, gives
\begin{equation}\label{nonvoid}
(\vec{n},\vec{k})\Omega^{ij}\,\nabla_{ij}\zeta_2+2(\nabla z,\nabla
\zeta_2)
=-\frac{\det(\nabla_{ij}\zeta_1)}{\det(\Omega_{ij})}-|\nabla
\zeta_1|^2,
\end{equation}
where $\zeta_1$ and $\zeta_2$ respectively is the third component of $\vec{\tau}_1$
and $\vec{\tau}_2$.
\begin{remark}\label{parziali}
\eqref{simply} is the linearization with respect to
\eqref{burstin}. In addition, \eqref{simply}, \eqref{nonvoid} and
the linearization with respect to \eqref{heinz} all are belong to
the same kind of degenerate elliptic equation.
\end{remark}

The degenerate elliptic equations we shall study is very closely
related to rigidity problems arising from infinitesimal isometric
deformation, as well as other geometry problem, such as minimal
surface in hyperbolic space, etc. In particular, the existence of
solutions with high order regularity is very important to
investigate many geometry problems. One would like to know under
what conditions the solution of such equations are as smooth as the
given data. The theory on well-posedness and regularity of solutions
to such equations, plays a crucial role in the above fields. Anyway,
such equations are deserved to be investigated vastly. However, so
far such equations might not be able to be treated by any standard
methods. Therefore maybe they will stimulate a general study of
linear, semilinear, quasilinear, and fully nonlinear degenerate
elliptic equations.

\subsection{The main question.}
The present paper is to devoted to investigate the well-posedness of
boundary value problems for a special class of degenerate linear
elliptic equations with previous geometric backgrounds. The aim of
this subsection is to bring up the main question of this paper. We
start with a few definitions and introduce notation and terminology
that is consistent throughout this paper.

Let $\Omega\subset \mathbb{R}^2$ be a bounded simply-connected
domain with smooth boundary $\partial\Omega$. And $\Omega$ is
divided into two subdomains by a smooth closed curve $\Gamma$. One
of the subdomains is called interior subdomain which is a
simply-connected one of $\Omega$. We denote it by $\Omega_+$, i.e.
$\Omega_+\subset\subset \Omega$. Another is denoted by
$\Omega_-:=\Omega\setminus\overline{\Omega}_+$, which is
connected. The boundary of $\Omega_+$ is denoted by $\partial
\Omega_+$. In addition, let $\varphi$ be a function in
$\overline{\Omega}$, and set
$$\Gamma:=\partial \Omega_+
=\{(\xi_1,\xi_2)\in \Omega\,\,\big|\,\,\varphi(\xi_1,\xi_2)=0\},$$
$$\Omega_+:=\{(\xi_1,\xi_2)\in\Omega\,\,\big|\,\,\varphi(\xi_1,\xi_2)>0\},$$
$$\Omega_-:=\Omega\setminus\overline{\Omega}_+=
\{(\xi_1,\xi_2)\in\Omega\,\,\big|\,\,\varphi(\xi_1,\xi_2)<0\},$$
where $\varphi$ is called the definition function of $\Gamma$.
Moreover, we suppose $\nabla\varphi\neq 0\quad
\hbox{on}\,\,\,\Gamma$, and denote the inward normal direction to
the boundary $\partial \Omega_+$ by $\vec{n}$. Obviously,
$$\vec{n}=
(n_1,n_2)= (\frac{\varphi_{\xi_1}}{|\nabla\varphi|},
\frac{\varphi_{\xi_2}}{|\nabla\varphi|})\big|_{\Gamma}\neq 0,$$
Consider
$$Lu\equiv\varphi (A^{ij}u_{\xi_i \xi_j}+C u)
+B^l u_{\xi_l}\quad \hbox{in}\,\,\Omega\subset\mathbb{R}^2_+,$$
where
\begin{equation}\label{fas}
\varphi,A^{ij},B^l,C\in C^\infty(\overline{\Omega})
\quad \hbox{for}\,\,\,i,j,l=1,2.
\end{equation}

Assume that
\begin{equation}\label{ken}
A^{ij} \eta_i \eta_j\geq\lambda_0|\eta|^2\,\quad
\hbox{for all}\,\,\,\eta=(\eta_1,\eta_2)\in \mathbb{R}^2,
\quad (\,i,j=1,2\,).
\end{equation}
where $\lambda_0$ is a positive constant;
\begin{eqnarray}\label{splin}
(B^l\varphi_{\xi_l})<0\quad\hbox{on}\,\,\,\Gamma;
\end{eqnarray}
and
\begin{equation}\label{patri}
C\leq 0\quad\hbox{for all}\,\,(\xi_1,\xi_2)\in \overline{\Omega}.
\end{equation}

Obviously, from the above assumptions it is easy to see that
\begin{equation}\label{iron}
(A^{ij}\varphi_{\xi_i} \varphi_{\xi_j})\big|_\Gamma
\geq\lambda_0(|\nabla\varphi|^2)\big|_\Gamma>0.
\end{equation}

Assume that $F\in L^2(\Omega)\,,\,\, g\in H^2(\Omega)\,.$ we shall
discuss the well-posedness of the following boundary value problem
(Abbreviation: BVP)
\begin{equation}\label{mesi}
\left\{
\begin{array}{lll}
Lu=F\quad \hbox{in}\,\,\,\Omega,\\
u=g\quad \hbox{on}\,\,\,\Gamma,\\
u=0\quad \hbox{on}\,\,\,\partial\Omega.
\end{array}
\right.
\end{equation}

Now we state the definition of weak solution in the following
\begin{definition}\label{egnaro}
Assume that $F\in L^2(\Omega), g\in H^2(\Omega)$.
$\mathcal{A}:H_0^1(\Omega)\times H_0^1(\Omega)\rightarrow
\mathbb{R}$ a continuous bilinear form, is defined by
$$\mathcal{A}(u,v)=\int_\Omega [-u_{\xi_i}(A^{ij}\varphi v)_{\xi_j}
+\varphi Cuv +B^lu_{\xi_l}v]d\xi \quad \hbox{for all}\,\,\,u,v\in
H^1_0(\Omega).$$ If $u\in H^1_0(\Omega)$ satisfies that
$$\left\{
\begin{array}{ll}
\mathcal{A}(u,v)=(F,v)\quad \hbox{for all}\,\,v\in H^1_0(\Omega),\\
u=g\quad \hbox{on}\,\,\, \Gamma,
\end{array}
\right.$$ where the boundary value is to be interpreted in the sense
of traces. Then $u$ is called $H^1$ weak solution of BVP
\eqref{mesi}.
\end{definition}

The question
\begin{quote}
\emph{Is there a $H^1$ solutions of BVP \eqref{mesi} and unique is
such solution?}
\end{quote}
is unknown as \emph{the well-posedness of the boundary value problem
for degenerate elliptic equations}.

Throughout this paper we will utilize such Convention:
\,\,(1)\,\,The C that are appearing in paper, all express positive
bounded constant. But they are possibly different when they are
appearing in different rows. \,\,(2)\,\,We often use
"$\rightharpoonup$" and "$\rightarrow$" expressing respectively
weak convergence and strong convergence in the corresponding
function spaces.

\subsection{The main result and its trivial generalization.}
The main result of this paper is the following
\begin{theorem}\label{lanoitanretni}
Suppose $\varphi$, $A^{ij}$, $B^l$\,\,\,$(i,j,l=1,2)$, and $C$
 satisfy the conditions \eqref{fas},
\eqref{ken}, \eqref{splin} and \eqref{patri}. Let $F\in L^2(\Omega)$
and $g\in H^2(\Omega)$. Then there exists an unique $H^1$ weak
solution $u$ of the BVP \eqref{mesi}, and $u$ satisfies
\begin{equation}\label{brod}
\|u\|_{H^1(\Omega)}\leq C\big[\,\|F\|_{L^2(\Omega)}
+\|g\|_{H^2(\Omega)}\,\big],
\end{equation}
where $C$ is a constant depending only on $\Gamma$,
$\|\varphi\|_{C^3(\overline{\Omega})}$,
$\|A^{ij}\|_{C^2(\overline{\Omega})}$,
$\|B^l\|_{C^1(\overline{\Omega})}$\quad $(i,j,l=1,2)$ and
$\|C\|_{C^1(\overline{\Omega})}$.
\end{theorem}
\begin{remark}\label{jenkins}
The BVP \eqref{mesi} can be discussed in $\mathbb{R}^{n+1}$ under
the same conditions. And the similar result also can be obtained by
the same methods in the case $\mathbb{R}^{n+1}$. Since all the
generalization is trivial, we omit to state this result and the
details of its proof at here.
\end{remark}
\begin{remark}
We explain briefly some known facts about the works of
Ole$\check{i}$nik and Radkevi$\check{c}$ (for details, to
see~\cite{OR}). Their theory requires that the characteristic form
of the equation is non-negative in the global domain. But the
characteristic form of equations in the problems which we deal with
is changing its signs in the domain. Next the theory of
Ole$\check{i}$nik and Radkevi$\check{c}$ requires that the
coefficient of unknown function term for the equation is negative
enough. It is usually not provided with this condition in the
practical problems. Hence we could not get the $L^2$ solution from
the direct applications of their conclusion to our problems. This is
the difficulties in our problems.
\end{remark}
Of course the regularity of solutions plays an important role in the
study of geometry problems. So we need to make a further discussion
on the corresponding regularity of solutions to such problems. The
further results on regularity will be given in our preprint
paper~\cite{HeYue}. In spite of many relevant progress, up to now
there has been no standard way to deal with such kind of problems,
and some crucial problems remain unsolved. Therefore, maybe our
methods are helpful in studying the general degenerate elliptic
equations.

\section{Preliminaries.}
\subsection{Homogenization of the BVP \eqref{mesi}.}
Suppose $g$ may be extend to domain $\overline{\Omega}$, still be
denoted by $g$, such that $g$ satisfies $g=0$\quad on\,\,\,$\partial
\Omega$. Without loss of generality we may assume that $g\equiv 0$.
In fact, if $u$ is a solution of BVP \eqref{mesi}, let $v=u-g$, then
$v$ is a solution of the following BVP:
$$\left\{
\begin{array}{lll}
Lv=F-Lg\quad \hbox{in}\,\,\,\Omega,\\
v=0\quad \hbox{on}\,\,\,\Gamma\cup\partial\Omega.
\end{array}
\right.$$  Hence, instead of the primary BVP \eqref{mesi} we may
discuss the well-posedness of the following BVP:
\begin{equation}\label{plenum}
\left\{
\begin{array}{ll}
  Lu=F\quad \hbox{in}\,\,\,\Omega,\\
  u=0\quad \hbox{on}\,\,\,\Gamma\cup\partial\Omega.
\end{array}
\right.
\end{equation}

\subsection{Simplification of the form to $Lu=F$.}
Firstly, we will simplify the form of the equation $Lu=F$ in some
neighborhood of $\Gamma$. By a appropriate transformations of the
variable in a neighborhood of $\Gamma$, we get the following result:
\begin{lemma}\label{albertcohen}
By a appropriate transformation $\Phi$ of the variable in some
neighborhood $\mathcal{N}_0(\Gamma)$ of $\Gamma$, $Lu=F$ can be
translated to the following form:
\begin{eqnarray}\label{nod}
\pounds u\equiv y(\omega u_{xx}+u_{yy}+cu)+au_x+bu_y=f, \quad
\hbox{on}\,\,\,D_{d_0},
\end{eqnarray}
where
$$D_{d_0}=\Phi\big(\mathcal{N}_0(\Gamma)\big)
=\{\,(x,y)\,\,\big|\,\,-\pi\leq x\leq \pi\,,\,-d_0<y<d_0\},$$ and
$\omega,a,b,c,f$ are all periodic functions with period $2\pi$ on
$x$. $\omega,a,b,c\in C^\infty(\bar{D}_{d_0})$.
\end{lemma}
\begin{proof}
We might as well suppose that $\Gamma$ may be expressed as follow
$$\Gamma=\{\,(\nu_1(s_1),\nu_2(s_1))\,\,\big|\,\,0\leq s_1\leq l\,\},$$
where $s_1$ is the parameter of arc length, and $l$ is the length
of $\Gamma$. Obviously,
$$\dot{\nu_1}^2+\dot{\nu_2}^2=1\quad \hbox{on}\,\,\,\Gamma,$$
where "$\cdot$" means the derivative with respect to $s_1$. Then
there exists a neighborhood $\mathcal{N}_0(\Gamma)$ of $\Gamma$,
such that any $\xi=(\xi_1,\xi_2)\in \mathcal{N}(\Gamma)$ may be
express as
\begin{equation}\label{transl}
\left\{
\begin{array}{ll}
  \xi_1=\nu_1(s_1)+n_1(s_1)s_2, \\
  \xi_2=\nu_2(s_1)+n_2(s_1)s_2,
\end{array}
\right.
\end{equation}
where $\vec{n}=(n_1,n_2)$ is the inward normal direction of
$\Omega_+$, Meanwhile it is also the exterior normal direction of
$\Omega_-$. In $\mathcal{N}_0(\Gamma)$, we know easily that
$\xi=(\xi_1,\xi_2)\in \Gamma$\,\,\, for\,\,\,$s_2=0$. On the other
hand, since $(\dot{\nu}_1,\dot{\nu}_2)$ is the unit tangent
direction of $\Gamma$, so $(\dot{\nu}_2,-\dot{\nu}_1)$ is the inward
normal direction of $\Gamma$. Hence
$(\dot{\nu}_2,-\dot{\nu}_1)=(n_1,n_2)$.

Clearly,
$$\aligned
\left.\frac{D(\xi_1,\xi_2)}{D(s_1,s_2)}\right|_{s_2=0}
&=\det\left.
\begin{pmatrix}
  \dot{\nu_1}+\dot{n}_1s_2 & n_1 \\
  \dot{\nu_2}+\dot{n}_2s_2 & n_2
\end{pmatrix}
\right|_{s_2=0}\\&=[\dot{\nu_1}n_2-\dot{\nu_2}n_1
+(\dot{\nu}_1\nu_2-\nu_1\dot{\nu}_2)s_2]\big|_{s_2=0} =n^2_1+n^2_2
=1.
\endaligned$$

Thus by the inverse function theorem, for $|s_2|\leq \delta$, there
exists a sufficient small constant $\delta>0$, such that $s_1,s_2$
are smooth functions on $\xi_1,\xi_2$. We express $s_1,s_2$ as
follow
$$\Xi:\quad
\left\{
\begin{array}{ll}
  s_1=s_1(\xi_1,\xi_2),\\
  s_2=s_2(\xi_1,\xi_2).
\end{array}
\right.$$
Applying \eqref{transl} again, we have
\begin{equation}\label{rosode}
\left\{
\begin{array}{ll}
  \displaystyle\frac{\partial s_1}{\partial \xi_1}=n_2/\bigtriangleup, \\[15pt]
  \displaystyle\frac{\partial s_2}{\partial \xi_1}=-(\dot{\nu}_2+n_2s_2)/\bigtriangleup,
\end{array}
\right. \quad\hbox{and}\quad \left\{
\begin{array}{ll}
  \displaystyle\frac{\partial s_1}{\partial \xi_2}=-n_1/\bigtriangleup, \\[15pt]
  \displaystyle\frac{\partial s_2}{\partial \xi_2}=(\dot{\nu}_1+n_1s_2)/\bigtriangleup,
\end{array}
\right.
\end{equation}
where $\bigtriangleup=\dot{\nu}_1n_2-\dot{\nu}_2n_1+(\dot{n}_1n_2-n_1\dot{n}_2)s_2$.

In the sequel, by the transform $\Xi$ and direct calculation, the
equation $Lu=F$ become
\begin{equation}\label{chinese}
\varphi\big(\widetilde{A}^{ij} u_{s_is_j}+Cu\big)
+\widetilde{B}^k u_{s_k}=F,
\end{equation}
where
\begin{equation}\label{etavired}
\widetilde{A}^{ij}
=A^{lr}\,\frac{\partial s_i}{\partial \xi_l} \frac{\partial s_j}{\partial \xi_r}\,,
\quad\quad
\widetilde{B}^k
=\varphi A^{ij}\,\frac{\partial^2 s_k}{\partial \xi_i \partial \xi_j}
+B^l \frac{\partial s_k}{\partial \xi_l}.
\end{equation}
Obviously, the transform
$\Xi:\mathcal{N}(\Gamma)\mapsto\Xi\big(\mathcal{N}(\Gamma)\big)$,
and makes $\Gamma$ to become
$$\Xi(\Gamma)=\{\,(s_1,s_2)\,\,\big|\,\,0\leq s_1\leq l\,,\,s_2=0\}.$$
So, we have
\begin{equation}\label{yhcuac}
\varphi(s_1,s_2)=s_2\widetilde{\varphi},
\end{equation}
where
$$\widetilde{\varphi}(s_1,s_2)=\int_0^1 \varphi_{s_2}(s_1,\varsigma
s_2)d\varsigma,$$

and
$$\aligned
\widetilde{\varphi}(s_1,0)&=\varphi_{s_2}(s_1,0)
\\[5pt]
&=\varphi_{\xi_1}\big|_\Gamma\,\frac{\partial \xi_1}{\partial
s_2}\big|_{s_2=0} +\varphi_{\xi_2}\big|_\Gamma\,\frac{\partial
\xi_2}{\partial s_2}\big|_{s_2=0}
\\[5pt]
&=(\varphi_{\xi_1}\,n_1+\varphi_{\xi_2}\,n_2)\big|_\Gamma\\[8pt]
&=\Big\{\sqrt{\varphi_{\xi_1}^2+\varphi_{\xi_2}^2}\,\,
\Big(\frac{\varphi_{\xi_1}}{\sqrt{\varphi_{\xi_1}^2+\varphi_{\xi_2}^2}}\,\,n_1
+\frac{\varphi_{\xi_2}}{\sqrt{\varphi_{\xi_1}^2+\varphi_{\xi_2}^2}}\,\,n_2\left.\Big)\Big\}
\right|_\Gamma
\\[10pt]
&=\left.\left(\sqrt{\varphi_{\xi_1}^2+\varphi_{\xi_2}^2}\,\,\right)\right|_\Gamma
\,\,(n_1^2+n_2^2)
\\[5pt]
&=(|\nabla\varphi|)\big|_\Gamma \neq 0.
\endaligned$$
Then there exists some neighborhood $U \subseteq
\Xi\big(\mathcal{N}(\Gamma)\big)$, such that
$\widetilde{\varphi}\neq 0$\quad in\,\,\,$U$.

Transforming the equation $Lu=F$ again, by
$$\Theta:\quad
\left\{\begin{array}{ll}
  x_1=\psi(s_1,s_2),\\
  x_2=s_2,
\end{array}\right.$$
where $\psi$ is a undetermined function. Then the form of the
equation $Lu=F$ is translated from \eqref{chinese} into the
following form:
$$x_2\,\widetilde{\varphi}\,\big[\big(\widetilde{A}^{ij}\, \displaystyle\frac{\partial x_p}{\partial s_i}
 \displaystyle\frac{\partial x_q}{\partial s_j}\big)u_{x_px_q}+C u\big]
+\big(\varphi \widetilde{A}^{ij}\,
 \displaystyle\frac{\partial^2 x_r}{\partial s_i \partial s_j}
 +\widetilde{B}^k \displaystyle\frac{\partial x_r}{\partial s_k}\big)u_{x_r}
=F.$$ By simplification, we obtain
\begin{eqnarray}\label{photo}
x_2\,\widetilde{\varphi}\,\big[(\widetilde{A}^{ij} \psi_{s_i} \psi_{s_j})u_{x_1x_1}
+(\widetilde{A}^{l2} \psi_{s_l})u_{x_1x_2}
+\widetilde{A}^{22}u_{x_2x_2}
+C u\big]&&\nonumber\\[5pt]
+(x_2 \widetilde{A}^{ij} \psi_{s_is_j}
+\widetilde{B}^k \psi_{s_k})u_{x_1}
+\widetilde{B}^2 u_{x_2}
&=&F.
\end{eqnarray}
In order to delete the mixed derivative term $(\widetilde{A}^{l2}
\psi_{s_l})u_{x_1x_2}$, we consider the Cauchy problem
\begin{equation}\label{munchen}
\left\{\begin{array}{ll}
\widetilde{A}^{12} \psi_{s_1}
+\widetilde{A}^{22} \psi_{s_2}=0,\\
 \psi(s_1,0)=s_1.
\end{array}
\right.
\end{equation}
Together \eqref{rosode} with \eqref{etavired}, we have
$$\widetilde{A}^{22}
  =(n_2,-n_1)
  \begin{pmatrix}
  A^{22} & -A^{12} \\
  -A^{12} & A^{11}
\end{pmatrix}
\begin{pmatrix}
  n_2 \\
  -n_1
\end{pmatrix}
>0,\quad \hbox{for}\,\,\,s_2=0.$$
Therefore,
$$(\widetilde{A}^{12},\widetilde{A}^{22})
\begin{pmatrix}
  0 \\
  1
\end{pmatrix}
=\widetilde{A}^{22}>0\quad \hbox{on}\,\,\,\Gamma.$$ On one hand,
according to the theory of one order PDE, there exists an unique
solution $\psi^*$ of the Cauchy problem \eqref{munchen} in a
neighborhood $V\,(\subseteq U)$ of $\{0\leq s_1\leq l,s_2=0\}$.
Obviously, $\psi^*(s_1,s_2)+l$ is a solution of the Cauchy problem
\begin{equation}\label{krantz}
\left\{\begin{array}{ll}
\widetilde{A}^{12} \psi_{s_1}
+\widetilde{A}^{22} \psi_{s_2}=0,\\
 \psi(s_1,0)=s_1+l.
\end{array}
\right.
\end{equation}
On the other hand, since $\widetilde{A}^{12}$ and
$\widetilde{A}^{22}$ are both periodic function with period $l$ on
$s$, thus $\psi^*(s_1+l,s_2)$ is also a solution of the Cauchy
problem \eqref{krantz}. Consequently, by the uniqueness of solution
of the Cauchy problem, $\psi^*(s_1+l,s_2)=\psi^*(s_1,s_2)+l$. By
choosing $\psi=\psi^*$, it follows that $\det\,J\big|_{s_2=0}=1>0$,
where
$$J=\begin{pmatrix}
  \psi^*_{s_1} & \psi^*_{s_2}\\
  0 & 1
\end{pmatrix}$$
is the Jacobi matrix of the transform $\Theta$. Additionally, by the
theory of one order PDE, we know that the solution depends on
continuously the initial data and the equation's coefficients. So
there exists a sufficient small neighborhood $W\,(\subseteq V)$ of
$\{0\leq s\leq l,s_2=0\}$, such that $\psi^*\in C^\infty(W)$ and
$\det\,J> 0\,\,\,\hbox{for all}\,\,(s_1,s_2)\in W$. Therefore, the
equation \eqref{photo} can be simplified as follow
$$y\,\widetilde{\varphi}\,[(\widetilde{A}^{ij} \psi^*_{s_i} \psi^*_{s_j})u_{x_1x_1}
+\widetilde{A}^{22}u_{x_2x_2} +C u]+(y \widetilde{A}^{ij}
\psi^*_{s_is_j} +\widetilde{B}^{k} \psi^*_{s_k})u_{x_1}
+\widetilde{B}^{2}u_{x_2}=F$$ in the neighborhood $W$ of $\{0\leq
s\leq l,s_2=0\}$.

Without loss of the generality, we may assume that
$$\Theta\circ\Xi(\Gamma)=\{\,(x_1,x_2)\,\,\big|\,\,-\pi\leq x_1\leq
\pi\,,\,x_2=0\},$$ and $$\Theta:W\mapsto\Theta(W) \subseteq \Theta
\circ \Xi\big(\mathcal{N}(\Gamma)\big) \subseteq
\{\,(x_1,x_2)\,\,\big| \,\,-\pi\leq x_1\leq \pi\}.$$ According to
all the above analyses, we know that $D_{d_0}\subseteq \Theta(W)$
for an adequate small $d_0$, where $$D_{d_0}=\{\,(x_1,x_2)\,\,\big|
\,\,-\pi\leq x_1\leq \pi\,,\,-d_0<x_2<d_0.\}.$$ Define
$\Phi=\Theta\circ\Xi$, and denote
$\mathcal{N}_0(\Gamma)=\Phi^{-1}(D_{d_0})$. Obviously,
$\mathcal{N}_0(\Gamma)\subseteq \mathcal{N}(\Gamma)$. Since
$\widetilde{\varphi}\big|_{\mathcal{N}_0(\Gamma)}>0,\,
\widetilde{A}^{22}\big|_{\mathcal{N}_0(\Gamma)}>0$, therefore
\eqref{photo} can be simplified as follow
\begin{eqnarray}\label{birkhauser}
\quad\quad \pounds u\equiv x_2(\omega
u_{x_1x_1}+u_{x_2x_2}+cu)+au_{x_1}+bu_{x_2}=f\quad \hbox{in}~~
D_{d_0},
\end{eqnarray}
where
\begin{equation}\label{narrate}
\left\{
\begin{array}{c}
\displaystyle \omega=\frac{\widetilde{A}^{ij}\psi^*_{s_i}
\psi^*_{s_j}} {\widetilde{A}^{22}},\,\,\,
a=\frac{x_2
\widetilde{A}^{ij} \psi^*_{s_is_j}+\widetilde{B}^k \psi^*_{s_k}}
{\widetilde{A}^{22}\,\widetilde{\varphi}},\,\,\,\\[15pt]
\displaystyle b=\frac{\widetilde{B}^{2}}
{\widetilde{A}^{22}\,\widetilde{\varphi}},\,\,\,
c=\frac{C}{\widetilde{A}^{22}},\,\,\,
f=\frac{F}{\widetilde{A}^{22}\,\,\widetilde{\varphi}}.
\end{array}
\right.
\end{equation}
For the simplicity, we still denote respectively $x_1,x_2$ as $x,y$
from now on. Thus we simply rewrite \eqref{birkhauser} as
\eqref{nod}. Obviously, from \eqref{narrate} it follows that
$\omega,a,b,c,f$ all are periodic function with period $2\pi$ on
$x$. In addition, from \eqref{fas} it follows that $\omega,a,b,c\in
C^\infty(\overline{D_{d_0}})$\,.
\end{proof}
For the sake of simplicity, denoting $\omega_0=\omega(x,0),\,\,
 a_{0}=a(x,0),\,\,
 b_{0}=b(x,0);$
$$L^\varepsilon u\equiv
(\varphi+\varepsilon) (A^{ij}u_{\xi_i\xi_j}+C u)+B^l u_{\xi_l}\quad (\,i,j,l=1,2\,);$$
$$L^{(-\varepsilon)} u\equiv
(\varphi-\varepsilon) (A^{ij}u_{\xi_i\xi_j}+C u)+B^l u_{\xi_l}\quad (\,i,j,l=1,2\,);$$
$$\pounds^\varepsilon u\equiv
(y+\varepsilon)(\omega u_{xx}+u_{yy}+cu)+au_x+bu_y;$$
$$\pounds^{(-\varepsilon)}u\equiv
(y-\varepsilon)(\omega u_{xx}+u_{yy}+cu)+au_x+bu_y;$$
$$D^+_d=\{\,(x,y)\,\,\big|\,\,-\pi\leq x\leq \pi\,,\,0<y<d\};$$
$$D^-_d=\{\,(x,y)\,\,\big|\,\,-\pi\leq x\leq \pi\,,\,-d<y<0\}.$$

In the sequel, we shall always use such convention: we sometimes
identify $u\circ\Phi^{-1}$ as $u$, and still denote
$u\circ\Phi^{-1}$ by $u$. Anyway, no confusion of ideas will rise in
this paper if only one keep concretely close touch with the context.
\begin{proposition}\label{actamath}
\begin{equation}\label{acinis}
b_0=\Big(\frac{B^l\varphi_{\xi_l}}
{A^{ij}\varphi_{\xi_i}\varphi_{\xi_j}}\Big)\Big|_\Gamma.
\end{equation}
\end{proposition}
\begin{proof}
The result follows from \eqref{etavired}, \eqref{narrate} and
\eqref{yhcuac}.
\end{proof}
\begin{remark}\label{jawo}
By \eqref{acinis} and \eqref{narrate}, it is easily verify that
$$b_0<0\quad\hbox{is equivalent to}
\quad\eqref{splin}\,\,\hbox{i.e.}\quad
(B^l\varphi_{\xi_l})<0\quad\hbox{on}\,\,\,\Gamma;$$
$$c\leq 0\quad\hbox{is equivalent to}
\quad\eqref{patri}\,\,\hbox{i.e.}\quad C\leq 0\quad\hbox{for
all}\,\,(\xi_1,\xi_2)\in \overline{\Omega}.$$
\end{remark}

\subsection{Elliptic regularization of the BVP \eqref{plenum}.}
Under the conditions \eqref{ken}, \eqref{splin} and \eqref{patri},
we will employ elliptic regularization to discuss the well-posedness
of the following BVP:
\begin{equation}\label{market}
\left\{
\begin{array}{ll}
Lu=F\quad \hbox{in}\,\,\,\Omega_+,\\
u=0\quad \hbox{on}\,\,\,\Gamma,
\end{array}
\right.
\end{equation}
and
\begin{equation}\label{radkeic}
\left\{
\begin{array}{lll}
Lu=F\quad \hbox{in}\,\,\,\Omega_-,\\
u=0\quad \hbox{on}\,\,\,\Gamma\cup\partial\Omega.
\end{array}
\right.
\end{equation}

Since $\varphi>0$,\,\,\,on\,\,\,$\Omega_+$;
$\varphi<0$\,\,\,on\,\,\,$\Omega_-$. then we thus may construct
the following subsidiary BVP:
\begin{equation}\label{scienze}
\left\{
\begin{array}{ll}
 L^\varepsilon u=F^\varepsilon\quad \hbox{in}\,\,\,\Omega_+,\\
 u=0\quad \hbox{on}\,\,\,\Gamma,
\end{array}
\right.
\end{equation}
and
\begin{equation}\label{oleinik}
\left\{
\begin{array}{lll}
L^{(-\varepsilon)} u=F^\varepsilon\quad \hbox{in}\,\,\,\Omega_-,\\
u=0\quad \hbox{on}\,\,\,\Gamma\cup\partial\Omega,
\end{array}
\right.
\end{equation}
where $F^{\varepsilon}\in \mathrm{C}^{\infty}(\overline{\Omega})$.
In fact, we may choose $F^{\varepsilon}$ as the mollification of $F$.
Moreover, if $F\in H^k(\Omega)$, then
\begin{equation}\label{bluee}
\|F^{\varepsilon}\|_{H^k(\Omega)}\leq C_k\|F\|_{H^k(\Omega)},
\end{equation}
where $C_k$ is a constant depending only on $k$; Furthermore,
$$\|F^{\varepsilon}-F\|_{H^k(\Omega)}\rightarrow 0,
\quad \hbox{as}\,\,\,\varepsilon\rightarrow 0. $$

\section{$H^1$ estimates for solutions of the BVP \eqref{scienze}.}
In this section, we will discuss the $H^1$ estimates for solutions
of the BVP \eqref{scienze}. According to the $L^2$ theory of second
order elliptic type equation, there exists a solution $u^\varepsilon
\in C^\infty(\overline{\Omega}_+)$ of the BVP \eqref{scienze}, and
also a solution $u^{(-\varepsilon)} \in
C^\infty(\overline{\Omega}_-)$ of the BVP \eqref{oleinik}.
Obviously, the interior $H^1$ estimates for solutions of above the
BVPs, can be derived from the standard interior $H^2$ estimates of
second order elliptic type equation. Therefore, we only need to give
the local estimates in a neighborhood of $\Gamma$. In addition, from
the Lemma \ref{albertcohen}, it follow that if we want to estimate
the solution $u^\varepsilon$\,\, (or $u^{(-\varepsilon)}$)\,\, of
the BVP \eqref{scienze}\,\,(or \eqref{oleinik})\,\, in a
neighborhood of $\Gamma$, then only need under the conditions
$b_0<0,\,c\leq 0$, to make the local $H^1$ estimates for the
solution of $\pounds^\varepsilon u=f^\varepsilon \quad
\hbox{in}\,\,\,D^+_{d_0}$,\,\,\,with\,\,\,$u(x,0)=0$\,\, (or
$\pounds^{(-\varepsilon)}u=f^\varepsilon \quad
\hbox{in}\,\,\,D^-_{d_0}$, \,\,\,with\,\,\,$u(x,0)=0$)\,\, in a
neighborhood of $\{y=0\}$, where $f^\varepsilon$ is the
mollification of $f$, satisfies that $\|f^{\varepsilon}\|_{L^2}\leq
C \|f\|_{L^2}$, and constant $C$ is independent of $\varepsilon$.

\subsection{$H^1$ estimates for the solution
of $\pounds^\varepsilon u=f^\varepsilon$ with $u(x,0)=0$} Under the
condition $b_0<0,\,c\leq 0$, we will give the $H^1$ estimates for
the solution $u^\varepsilon$ of $\pounds^\varepsilon
u=f^\varepsilon$ with $u(x,0)=0$ in the following
\begin{proposition}\label{film}
Suppose $\omega\geq \omega_*>0\,; b_0<0,\,c\leq 0$. Then there
exists a constant $\sigma=\sigma(\omega,a,b)\in(0,d_0/2]$, such that
for any $d\in(0,\sigma]$, the solution $u^\varepsilon$ of the
equation $\pounds^\varepsilon u=f^\varepsilon$ with $u(x,0)=0$ has
the following estimates
\begin{equation}\label{muggy}
\|\nabla u^\varepsilon\|_{L^2(D^+_d)}\leq C
\big[\,\|f^\varepsilon\|_{L^2(D^+_{2d})}
+\|u^\varepsilon\|_{L^2(D^+_{2d})}\,\big],
\end{equation}
where $C$ is a constant depending only on $\omega$, $a$, $b$ and $d$.
\end{proposition}
\begin{proof}
First, we construct the cutoff function $\vartheta\in
C^\infty(\mathbb{R}_+)$, such that $0\leq\vartheta\leq 1,$ and
satisfies
\begin{equation}\label{wine}
\vartheta=\vartheta(y)=\left\{\begin{array}{ll}
    1\,,\,\,0\leq y\leq d,\\
    0\,,\,\,y\geq 2d,
  \end{array}\right.
  \quad \hbox{and}\,
  \quad \big|\frac{\partial^k\vartheta}{\partial y^k}\big|\leq\frac{C}{d^k}\,,
  \quad k\in \mathbb{N}.
  \end{equation}
Define $v^\varepsilon=\vartheta u^\varepsilon.$ Thus, from
$\pounds^\varepsilon u^\varepsilon=f^\varepsilon$ it follows
\begin{equation}\label{you}
\pounds^\varepsilon v^\varepsilon=
\vartheta f^\varepsilon
+2(y+\varepsilon)\vartheta_yu^\varepsilon_y
+[(y+\varepsilon)\vartheta_{yy}+b\vartheta_y]u^\varepsilon
\quad \hbox{in}\,\,\,D^+_{d_0}.
\end{equation}
Clearly, $u^\varepsilon(x,0)=0$ implies $v^\varepsilon(x,0)=0$.

Now we make both sides of \eqref{you} inner product with
$v^\varepsilon/(y+\varepsilon)$, i.e.,
\begin{equation}\label{theory}
-\big(\pounds^\varepsilon
v^\varepsilon,\frac{v^\varepsilon}{y+\varepsilon}\big)
=-\big(\vartheta
f^\varepsilon+2(y+\varepsilon)\vartheta_yu^\varepsilon_y
+[(y+\varepsilon)\vartheta_{yy} +b\vartheta_y]u^\varepsilon
,\frac{v^\varepsilon}{y+\varepsilon}\big).
\end{equation}
Moreover, denoting
$$\iint:=\int_{D^+_d}:=\int^\pi_{-\pi}\int^d_0.$$
For the left-hand side of $\eqref{theory}$, we get
$$\aligned
&-\big(\pounds^\varepsilon
v^\varepsilon,\frac{v^\varepsilon}{y+\varepsilon}\big)\\
=&\iint\big[-(\omega
v^\varepsilon_{xx}+v^\varepsilon_{yy})v^\varepsilon
-c|v^\varepsilon|^2
-\frac{a}{y+\varepsilon}\,v^\varepsilon_xv^\varepsilon
-\frac{b}{y+\varepsilon}\,v^\varepsilon_yv^\varepsilon\big]\\[10pt]
=&\iint\big[-(\omega v^\varepsilon_xv^\varepsilon)_x
-(v^\varepsilon_yv^\varepsilon)_y -c|v^\varepsilon|^2
+\omega_xv^\varepsilon_xv^\varepsilon +\omega |v^\varepsilon_x|^2
+|v^\varepsilon_y|^2\nonumber\\[10pt]&
+\frac{a_x}{2y+2\varepsilon}|v^\varepsilon|^2
+\big(\frac{b}{2y+2\varepsilon}\big)_y|v^\varepsilon|^2\big]\\[10pt]&
-\int^{2d}_0\frac{1}{2y+2\varepsilon}\,(a|v^\varepsilon|^2)\big|^{2d}_{-2d}dy
-\int^{2d}_{-2d}\big(\frac{b|v^\varepsilon|^2}{2y+2\varepsilon}\big)\big|^{2d}_0dx
\\[10pt]
=&\iint\Big\{\omega |v^\varepsilon_x|^2 +|v^\varepsilon_y|^2
+\omega_xv^\varepsilon_xv^\varepsilon
+\big[\frac{a_x}{2y+2\varepsilon}
+\big(\frac{b}{2y+2\varepsilon}\big)_y-c\big]|v^\varepsilon|^2\Big\}
\\[10pt]
=&\iint\Big\{|v^\varepsilon_x|^2 +|v^\varepsilon_y|^2
+\omega_xv^\varepsilon_xv^\varepsilon +\big[\frac{(y+\varepsilon)(
a_x+b_y)-b} {2(y+\varepsilon)^2}-c\big]|v^\varepsilon|^2\Big\}.
\endaligned$$
and
$$\iint \omega_xv^\varepsilon_xv^\varepsilon\leq
\frac{1}{4}\iint \omega |v^\varepsilon_x|^2
+\iint\vartheta^2\frac{\omega_x^2}{\omega}|u^\varepsilon|^2.$$ In
addition, the terms in the right-hand side of \eqref{theory} have
the following estimates:
$$\iint\vartheta f^\varepsilon \frac{v^\varepsilon}{y+\varepsilon}
\leq\frac{4}{\delta}\iint \vartheta^2|f^\varepsilon|^2
+\frac{\delta}{16}\iint\frac{|v^\varepsilon|^2}{(y+\varepsilon)^2};$$
\vspace{5pt}
$$\aligned
\iint
2(y+\varepsilon)\,\vartheta_yu^\varepsilon_y\,\frac{v^\varepsilon}{y+\varepsilon}
=& \iint 2\vartheta_y(\vartheta u^\varepsilon_y)u^\varepsilon
\nonumber\\[5pt]
=&\iint 2\vartheta_y[(\vartheta
u^\varepsilon)_y-\vartheta_yu^\varepsilon]u^\varepsilon
\nonumber\\[5pt]
=&\iint
2(\vartheta_yv^\varepsilon_yu^\varepsilon-\vartheta_y^2|u^\varepsilon|^2),
\endaligned$$
and
$$\iint \vartheta_yv^\varepsilon_yu^\varepsilon\leq \iint
\vartheta_y^2|u^\varepsilon|^2+\frac{1}{4}\iint
|v^\varepsilon_y|^2,$$
Thus we have
$$\iint 2(y+\varepsilon)
\,\vartheta_yu^\varepsilon_y\,\frac{v^\varepsilon}{y+\varepsilon}
\leq \frac{1}{2}\iint |v^\varepsilon_y|^2;$$
$$\iint
b\vartheta_y\,u^\varepsilon\frac{v^\varepsilon}{y+\varepsilon}
\leq \frac{4}{\delta}\iint
|b|_\infty^2\vartheta_y^2|u^\varepsilon|^2
+\frac{\delta}{16}\iint\frac{|v^\varepsilon|^2}{(y+\varepsilon)^2};$$
and
$$\iint (y+\varepsilon) \vartheta_{yy}u^\varepsilon\frac{v^\varepsilon}{y+\varepsilon}
\leq \iint \vartheta|\vartheta_{yy}||u^\varepsilon|^2.$$ Together
with \eqref{theory} and  all above estimates, it implies
$$\aligned
&\iint\Big\{\omega |v^\varepsilon_x|^2+|v^\varepsilon_y|^2
+\big[\frac{(y+\varepsilon)( a_x+b_y)-b}
{2(y+\varepsilon)^2}-c\big]|v^\varepsilon|^2\Big\}
\\[10pt]\leq&
\frac{4}{\delta}\iint \vartheta^2|f^\varepsilon|^2
+\iint \big(\frac{4}{\delta} |b|_\infty^2 \vartheta_y^2
+\vartheta^2\frac{\omega_x^2}{\omega}
+\vartheta |\vartheta_{yy}| \big)|u^\varepsilon|^2\\[10pt]&
+\frac{1}{4}\iint \omega |v^\varepsilon_x|^2 +\frac{1}{2}\iint
|v^\varepsilon_y|^2
+\frac{\delta}{8}\iint\frac{|v^\varepsilon|^2}{(y+\varepsilon)^2}.
\endaligned$$
Hence by a simplification procedure we obtain
\begin{eqnarray}\label{newyork}
&&\iint
\Big\{\frac{3\omega}{2}|v^\varepsilon_x|^2+|v^\varepsilon_y|^2
+\big[\frac{(y+\varepsilon)( a_x+b_y)-b-\frac{\delta}{4}}
{2(y+\varepsilon)^2}-c\big]|v^\varepsilon|^2\Big\}
\\[10pt]&\leq&
\frac{8}{\delta}\iint \vartheta^2|f^\varepsilon|^2 +\iint
\big(\frac{8}{\delta} |b|_\infty^2 \vartheta_y^2
+2\vartheta^2\frac{\omega_x^2}{\omega} +2\vartheta |\vartheta_{yy}|
\big)\,|u^\varepsilon|^2.\nonumber
\end{eqnarray}
Since $b_0<0$. By the continuity of $b_0$, there exists a constant
$\delta_0>0$ such that $-b_0> \delta_0$. Subsequently, by the
continuity of $b$, there exists an adequate small constant
$\sigma_1$ such that $-b\geq \frac{\delta_0}{2}$ for $0\leq
y\leq\sigma_1$. So, we have
$-b-\frac{\delta_0}{4}>\frac{\delta_0}{4}$ for $0\leq
y\leq\sigma_1$. Choose $\delta=\delta_0$, from \eqref{newyork} we
get
\begin{eqnarray}\label{girl}
&&\iint \Big\{\frac{3\omega}{2}|v^\varepsilon_x|^2+|v^\varepsilon_y|^2
+\big[\frac{(y+\varepsilon)
(a_x+b_y)+\frac{\delta_0}{4}}{2(y+\varepsilon)^2}-c\big]|v^\varepsilon|^2\Big\}
\\[10pt]&\leq&
\frac{8}{\delta_0}\iint \vartheta^2|f^\varepsilon|^2 +\iint
\big(\frac{8}{\delta_0} |b|_\infty^2 \vartheta_y^2
+2\vartheta^2\frac{\omega_x^2}{\omega} +2\vartheta
|\vartheta_{yy}|\big)\,|u^\varepsilon|^2.\nonumber
\end{eqnarray}
Obviously, we may choose a sufficient small constant $\sigma_2>0$,
which only depends on $\delta_0,a,b$; such that
$(y+\varepsilon)(a_x+b_y)+\frac{\delta_0}{4}\geq 0$, for $0\leq
y\leq \sigma_2$.

Choose a constant $\sigma=\min\{\sigma_1,\sigma_2,d_0/2\}$ and
together $c\leq 0$ with \eqref{girl}, we therefore obtain that
$$\iint \big(\frac{3\omega}{2}|v^\varepsilon_x|^2+|v^\varepsilon_y|^2\big)
\leq \frac{8}{\delta_0}\iint \vartheta^2|f^\varepsilon|^2 +\iint
\big(\frac{8}{\delta_0} |b|_\infty^2 \vartheta_y^2
+2\vartheta^2\frac{\omega_x^2}{\omega} +2\vartheta |\vartheta_{yy}|
\big)\,|u^\varepsilon|^2.$$ always holds for any $d\leq \sigma$.
Finally, by utilizing $\omega\geq \omega_*>0$ and \eqref{wine}, we
have
\begin{eqnarray}\label{baby}
\int_{D^+_d}|\nabla u^\varepsilon|^2\leq C(\omega,b) \int_{D^+_{2d}}
|f^\varepsilon|^2 +\frac{C(\omega,a,b)}{d^2} \int_{D^+_{2d}}
|u^\varepsilon|^2.\nonumber
\end{eqnarray}
This implies \eqref{muggy}.
\end{proof}
\begin{proposition}[Boundary $H^1$ estimates]\label{skaut}
Assume the coefficients $\varphi,A^{ij},B^l$ $(i,j,l=1,2.)$ and $C$
of the operator $L$ satisfy the conditions
\eqref{fas},\,\eqref{ken},\,\eqref{splin} and \eqref{patri}. Then
the solution $u^\varepsilon$ of the BVP \eqref{scienze}, satisfies
that
\begin{equation}\label{quell}
\|u^\varepsilon\|_{H^1(\Phi^{-1}(D^+_\sigma))}
\leq C\big[\,\|u^\varepsilon\|_{L^2(\Omega_+)}+\|F^\varepsilon\|_{L^2(\Omega_+)}\,\big],
\end{equation}
where $C$ is a constant depending only on
$\|\varphi\|_{C^3(\overline{\Omega}_+)},\|A^{ij}\|_{C^2(\overline{\Omega}_+)}$
and $\|B^l\|_{C^1(\overline{\Omega}_+)}$ $(i,j,l=1,2)$.
\end{proposition}
\begin{proof}
Rewriting the inequality \eqref{muggy}
in term of the variable $\xi=(\xi_1,\xi_2)$, we have
$$\|\nabla u^\varepsilon\|_{L^2[\Phi^{-1}(D^+_\sigma)]}\leq
C \big\{\|u^\varepsilon\|_{L^2[\Phi^{-1}(D^+_{2\sigma})]}
+\|F^\varepsilon\|_{L^2[\Phi^{-1}(D^+_{2\sigma})]}\big\}. $$ Let
the integral domain $\Phi^{-1}(D^+_{2\sigma})$ of the right-hand
side of the above inequality extend to $\Omega_+$. The proof is
complete.
\end{proof}
\begin{lemma}[Interior $H^1$ estimates]\label{wang}
Under the conditions of Proposition \ref{skaut}. Then for arbitrary
$\Omega'\subset\subset\Omega_+$, the solution $u^\varepsilon$ of the
equation $L^\varepsilon u=F^\varepsilon$ has the following interior
estimates
\begin{equation}\label{shi}
\|u^\varepsilon\|_{H^1(\Omega')}\leq
C\big[\,\|u^\varepsilon\|_{L^2(\Omega_+)}+\|F^\varepsilon\|_{L^2(\Omega_+)}\,\big],
\end{equation}
where $C$ is a constant depending only on $dist\{\Omega',\Gamma\},
\|\varphi\|_{C^2(\overline{\Omega}_+)},
\|A^{ij}\|_{C^1(\overline{\Omega}_+)}$,
$\|B^l\|_{C^1(\overline{\Omega}_+)}$ and
$\|C\|_{C^1(\overline{\Omega}_+)}$.
\end{lemma}
\begin{proof}
From the interior $H^2$ estimates of second order elliptic type
equation directly follows the conclusion.
\end{proof}
\begin{proposition}[Global $H^1$ estimates]\label{skeil}
Under the conditions of Proposition \ref{skaut}. Then the solution
$u^\varepsilon$ of BVP \eqref{scienze}, has the following global
estimates
\begin{equation}\label{maiti}
\|u^\varepsilon\|_{H^1(\Omega_+)}
\leq C\big[\,\|u^\varepsilon\|_{L^2(\Omega_+)}+\|F^\varepsilon\|_{L^2(\Omega_+)}\,\big],
\end{equation}
where $C$ is a constant depending only on
$\|\varphi\|_{C^3(\overline{\Omega}_+)},\|A^{ij}\|_{C^2(\overline{\Omega}_+)}$,
$\|B^l\|_{C^1(\overline{\Omega}_+)}$ $(i,j,l=1,2)$ and
$\|C\|_{C^1(\overline{\Omega}_+)}$.
\end{proposition}
\begin{proof}
Choose $\Omega'=\Omega_+ \setminus \Phi^{-1}(D^+_{\sigma/2})$.
\eqref{maiti} follows from \eqref{quell} and \eqref{shi}.
\end{proof}

\subsection{$H^1$ estimates of solutions of the BVP \eqref{market}.}
\begin{lemma}[see Lemma 3.1 of \cite{H}]\label{heyue}
Let $f\in H^1(D_{d_0}^+)$, if $v\in H^1(D_{d_0}^+)$ is a weak
solution of the following BVP:
$$\left\{
\begin{array}{lll}(L-\lambda I)v=f\quad \hbox{in}\,\,\,D_{d_0}^+,\\
v\,\,\hbox{is a periodic function with period}\,\,2\pi\,\,\hbox{on}\,\,x,\\
v(x,0)=v(x,d_0)=0.
\end{array}\right.$$
Then there exists a $\lambda_0>0$, such that when $\lambda\geq
\lambda_0$, we have $v_x\in H^1(D_{d_0}^+)$, and
$$\|v_x\|_1\leq C \|f\|_1,$$
where $C$ is a constant.
\end{lemma}
\begin{lemma}\label{qiyueershijiu}
Let $u,u_x\in H^1(\Omega)$. Then $u\in C(\overline{\Omega})$.
\end{lemma}
\begin{proof}
In fact, by the localization technique, we only need prove one conclusion:

let $u,u_x\in H^1(\mathbb{R}^2)$, then
$u\in C(\mathbb{R}^2)$.

Firstly, that $u,u_x\in H^1(\mathbb{R}^2)$ is equivalent to
$$\frac{1}{2\pi}\int_{\mathbb{R}^2} |\widehat{u}|^2
[1+\xi^2+\eta^2+\xi^2(\xi^2+\eta^2)]d\xi d\eta<+\infty.$$
where $\widehat{u}$ expresses the Fourier transform of $u$.

Secondly, it is easy to verify
$$\frac{1}{2\pi}\int_{\mathbb{R}^2}
\frac{1}{1+\xi^2+\eta^2+\xi^2(\xi^2+\eta^2)}d\xi d\eta
<+\infty.$$
So we have
$$\aligned
&\frac{1}{2\pi}\int_{\mathbb{R}^2} |\widehat{u}|d\xi d\eta\\
=&\frac{1}{2\pi}\int_{\mathbb{R}^2}
|\widehat{u}|\,[1+\xi^2+\eta^2+\xi^2(\xi^2+\eta^2)]^\frac{1}{2}
\frac{1}{[1+\xi^2+\eta^2+\xi^2(\xi^2+\eta^2)]^\frac{1}{2}}d\xi d\eta\\
\leq&\Big(\frac{1}{2\pi}\int_{\mathbb{R}^2}
|\widehat{u}|^2[1+\xi^2+\eta^2+\xi^2(\xi^2+\eta^2)]d\xi
d\eta\Big)^\frac{1}{2}\\
&\times \Big(\frac{1}{2\pi}\int_{\mathbb{R}^2}
\frac{1}{1+\xi^2+\eta^2+\xi^2(\xi^2+\eta^2)}d\xi
d\eta\Big)^\frac{1}{2}
<+\infty.
\endaligned$$
Thirdly, since
$u(x,y)=\frac{1}{2\pi}\int_{\mathbb{R}^2}
e^{i(x\xi+y\eta)}\widehat{u}(\xi,\eta)d\xi d\eta$,
thus we have
$$|u(x,y)-u(x_0,y_0)|
\leq \frac{1}{2\pi}\int_{\mathbb{R}^2}
|e^{i(x\xi+y\eta)}-e^{i(x_0\xi+y_0\eta)}|\cdot
|\widehat{u}(\xi,\eta)|d\xi d\eta.$$
Finally, we define $G(\xi,\eta;x,y;x_0,y_0)
=|e^{i(x\xi+y\eta)}-e^{i(x_0\xi+y_0\eta)}|\cdot|\widehat{u}(\xi,\eta)|$.
Obviously,
 $$G(\xi,\eta;x,y;x_0,y_0)\rightarrow 0\quad \hbox{a.e. in}\,\,\,\mathbb{R}^2
\quad \hbox{as}\,\,\,(x,y)\rightarrow (x_0,y_0)$$
In addition, $|G(\xi,\eta;x,y;x_0,y_0)|\leq 2|\widehat{u}(\xi,\eta)|
\in L^1(\mathbb{R}^2)$.
Hence the Lebesgue Dominated Convergence Theorem implies
$$|u(x,y)-u(x_0,y_0)|\rightarrow 0\quad \hbox{as}\,\,\,(x,y)\rightarrow (x_0,y_0).$$
i,e., $u\in C(\mathbb{R}^2)$.
\end{proof}
\begin{lemma}\label{weiyixing}
Let $u^*$ be a $H^1$ weak solution of the following BVP:
\begin{equation}\label{noitazilacol}
\left\{
\begin{array}{ll}
  Lu=0\quad \hbox{in}\quad \Omega_+,\\
  u=0\quad \hbox{on}\quad \Gamma.
\end{array}
\right.
\end{equation}
Then $u^*\equiv 0$.
\end{lemma}
\begin{proof}
Firstly, we prove $u^*\in C(\overline{\Omega}_+)\cap
C^\infty(\Omega_+)$. According to the interior regularity of
second order elliptic type equation, we know that $u^*\in
C^\infty(\overline{\Omega'}),\,\, \hbox{for
all}\,\,\Omega'\subset\subset\Omega_+$. Because of $L$ only
degenerates on $\Gamma$. Thus we only discuss the continuity of
$u^*$ in some neighborhood of $\Gamma$.

From the Lemma \ref{albertcohen}, it is easy to know that $Lu^*=0$
may be simplified as $\pounds u^*=0$ by the transform
$\Phi:\mathcal{N}_0(\Gamma)\mapsto D_{d_0}$. Define a cutoff
function $\zeta\in C^\infty(\mathbb{R}^1_+)$ as follows
$$\zeta=\zeta(y)=
\left\{
\begin{array}{ll}
  1\quad \hbox{for}\,\,\,0\leq y\leq \delta,\\
  0\quad \hbox{for}\,\,\,y\geq 2\delta.
\end{array}
\right.
\quad\quad \big(\,\hbox{for}\,\,\,\delta<\frac{d_0}{2}\,\big)$$
From $\pounds u^*=0$, it follows easily that
$(\pounds-\lambda)(\zeta u^*)=f^*$,
where $f^*=2y\zeta_yu^*_y+(y\zeta_{yy}+b\zeta_y-\lambda\zeta)u^*$.
Obviously, $\zeta u^*$ is a $H^1$ weak solution of problem
$$\left\{
\begin{array}{lll}
  (\pounds-\lambda)w=f^*\quad\hbox{in}\,\,\,D^+_{d_0}\\
  w\,\,\hbox{is a periodic function with period}\,\,2\pi\,\,\hbox{on}\,\,x,\\
  w(x,0)=w(x,d_0)=0.
\end{array}
\right.$$ Obviously, $u^*\in H^2(D^+_{d_0}\setminus
\overline{D^+_\delta})$. Direct calculation shows that $f^*\in
H^1(D^+_{d_0})$. Hence by applying Lemma \ref{heyue}, we have
$(\zeta u^*)_x\in H^1(D^+_{d_0})$. Thus Lemma \ref{qiyueershijiu}
implies $\zeta u^*\in C(\overline{D^+_{d_0}})$. Consequently,
$u^*\in C(\overline{D^+_\delta})$. Set
$\mathcal{N}^+_1(\Gamma)=\Phi^{-1}(D^+_\delta)$. Then $u^*\in
C\big(\overline{\mathcal{N}^+_1(\Gamma)}\big)$. Using the interior
regularity of elliptic equation, we deduce $u^*\in
C(\overline{\Omega}_+)\cap C^\infty(\Omega_+)$.

Next we will prove $u^*\equiv 0$. Since $u^*$ on $\Gamma=0$ and
$u^*\in C(\overline{\Omega}_+)$. for any $\epsilon>0$, there exists
a constant $\delta_0=\delta_0(\varepsilon)$ such that
$$|u^*(\xi,\eta)|<\epsilon\quad \hbox{for}\,\,\,dist\big((\xi,\eta),\Gamma\big)<\delta_0.$$
Choose simple connected domain $\Omega'\subset\subset\Omega_+$, such
that $dist(\partial\Omega',\Gamma)<\frac{\delta_0}{2}$. By the
maximum principle of elliptic equation, then $u^*$ must attains its
maximum and minimum of  $\overline{\Omega'}$ on $\partial\Omega'$,
i.e., $\max_{\overline{\Omega'}}|u^*|\leq
\max_{\partial\Omega'}|u^*|.$

Therefore,
$$\aligned
\max_{\overline{\Omega}_+} |u^*|\leq
\max\big\{\max_{\overline{\Omega'}}|u^*|,
\max_{\overline{\Omega\setminus\Omega'}}|u^*|\big\}
\leq \max\big\{\max_{\partial\Omega'}|u^*|,
\max_{\overline{\Omega\setminus\Omega'}}|u^*|\big\}
<\epsilon.
\endaligned$$
By the arbitrariness of $\epsilon$, we deduce that
$\max_{\overline{\Omega}_+} |u^*|\leq 0$.
This implies $u^*\equiv 0$.
\end{proof}
Now we  consider a family of elliptic operators
\[\Psi=\{\,L^\varepsilon\,\,\big|\,\, L^\varepsilon=
(\varphi+\varepsilon)(A^{ij}\,\frac{\partial^2}{\partial \xi_i
\partial\xi_j}+C) +B^l\,\frac{\partial}{\partial \xi_l};\,\,
\varphi,A^{ij},B^l,\,\,\hbox{and}\,\,C\,\,\]
\[\hbox{satisfy}\,\,\eqref{fas},\,\eqref{ken},
\,\eqref{patri};\,\,i,j,l=1,2,\,\,0<\varepsilon\leq 1.\}.\]
\begin{lemma}\label{pinl}
Suppose that for any $L^\varepsilon\in \Psi,\Gamma\in C^{1,1}$. Let
$u\in H^1(\Omega_+)$ with $u=0$ on ${\Gamma}$, and satisfies the
estimates
\begin{equation}\label{rei}
\|u\|_{H^1(\Omega_+)}
\leq C\big[\,\|u\|_{L^2(\Omega_+)}+\|L^\varepsilon u\|_{L^2(\Omega_+)}\,\big],
\end{equation}
Then we have
\begin{equation}\label{tardy}
\|u\|_{H^1(\Omega_+)}
\leq C\|L^\varepsilon u\|_{L^2(\Omega_+)},
\end{equation}
where constant $C$ is independent of $u$ and $\varepsilon$.
\end{lemma}
\begin{proof}
If the conclusion does not hold, then for any $n\in \mathbb{N}$,
there exist sequences $\{u_n\}$ satisfying
the assumed conditions and $\{\varepsilon_n\}$, such that
\begin{equation}\label{swift}
\|u_n\|_{H^1}\geq n\,\|L^{\varepsilon_n} u_n\|_{L^2}
\quad \hbox{and}\,\,\,\|u_n\|_{L^2}=1.
\end{equation}
By \eqref{rei}, it follows that
$$\|u_n\|_{H^1}\leq C\big(\frac{\|u_n\|_{H^1}}{n}+1\big).$$
Thus we have
\begin{equation}\label{straid}
\|u_n\|_{H^1}\leq C\quad \hbox{for}\,\,\,n\geq 2C.
\end{equation}
Therefore a subsequence $\{u_{n_k}\}$ of $\{u_n\}$ converges weakly in $H^1(\Omega_+)$
to some function $u^*\in H^1(\Omega_+)$, i.e.,
$$u_{n_k}\rightharpoonup u^*\quad \hbox{weakly in}\,\,\,H^1(\Omega_+)
\quad \hbox{as}\,\,\,k\rightarrow \infty.$$ According to the
Banach-Saks Theorem (cf. \cite{Ri}), we may choose a subsequence
of $\{u_{n_k}\}$, might as well still denote by $\{u_{n_k}\}$,
such that $\{\tilde{u}_k\}$ which is consist of the arithmetic
mean of $\{u_{n_k}\}$ as follow
$$\tilde{u}_k=\frac{1}{k}(u_{n_1}+\cdots+u_{n_k})$$
converges strongly in $H^1(\Omega_+)$ to $u^*\in H^1(\Omega_+)$,
i.e.,
$$\tilde{u}_k\rightharpoonup u^*\quad \hbox{strongly
in}\,\,\,H^1(\Omega_+)\quad \hbox{as}\,\,\,k\rightarrow \infty.$$
Obviously,
$$\tilde{u}_k=\frac{1}{k}(u_{n_1}+\cdots+u_{n_k})=0\quad \hbox{on}\,\,\,\Gamma.$$
By the Trace Theorem, it follows that
$u^*=0\quad \hbox{on}\,\,\,\Gamma$, in the sense of traces.

Since $\{\varepsilon_{n_k}\}$ is bounded, there exists a
subsequence of $\{\varepsilon_{n_k}\}$,
might as well still denote by $\{\varepsilon_{n_k}\}$, such that
$$\varepsilon_{n_k}\rightarrow \varepsilon_0\quad \hbox{as}\,\,\,k\rightarrow
\infty.$$

In addition, by the Sobolev embedding theorem,
$$u_{n_k}\rightarrow u^*\quad \hbox{strongly in}\,\,\,L^2(\Omega_+)
\quad \hbox{as}\,\,\,k\rightarrow \infty.$$ So we have
$$L^{\varepsilon_{n_k}} u_{n_k}\rightarrow L^{\varepsilon_0} u^*
\quad \hbox{as}\,\,\,k\rightarrow \infty.$$ in the sense of
distribution. Together \eqref{swift} with \eqref{straid}, we deduce
$$\|L^{\varepsilon_{n_k}} u_{n_k}\|_{L^2}\leq \frac{C}{n_k},$$
Letting $k\rightarrow \infty$, we conclude $\|L^{\varepsilon_0}
u^*\|_{L^2}=0$. Hence in the sense of distribution, we obtain
$L^{\varepsilon_0} u^*=0$.

Since $u^*=0\quad \hbox{on}\,\,\,\Gamma$ and $L^{\varepsilon_0}
u^*=0$, thus $u^*$ is a $H^1$ weak solution of the following BVP:
\begin{equation}\label{iveid}
\left\{
\begin{array}{ll}
  L^{\varepsilon_0} u=0\quad \hbox{in}\,\,\,\Omega_+,\\
  u=0,\quad \hbox{on}\,\,\,\Gamma.
\end{array}
\right.
\end{equation}
Then we have $u^*\equiv 0$.

We will respectively prove $u^*\equiv 0$ in the following two cases:

1)\,\,$\varepsilon_0>0$, then by the Maximal Principle of second
order elliptic type equation, we know that the BVP \eqref{iveid}
only has null solution. So, $u^*=0$.

2)\,\,$\varepsilon_0=0$, i.e., $u^*$ is a $H^1$ weak solution of
the BVP \eqref{noitazilacol}. Then  Lemma \ref{weiyixing} implies
$u^*\equiv 0$.

On the other hand,
$$1=\|u_{n_k}\|_{L^2}\rightarrow \|u^*\|_{L^2}\quad \hbox{as}\,\,\,k\rightarrow\infty,$$
yields $\|u^*\|_{L^2}=1$.
This contradicts $u^*=0$.
The proof for \eqref{tardy} is complete.
\end{proof}
\begin{theorem}[Global $H^1$ estimates]\label{ritein}
Under the conditions of Proposition \ref{skaut}. Then the solution
$u^\varepsilon$ of the BVP \eqref{scienze}, has the following global
estimates
\begin{equation}\label{fikst}
\|u^\varepsilon\|_{H^1(\Omega_+)}
\leq C\|F\|_{L^2(\Omega_+)},
\end{equation}
where $C$ is a constant depending only on
$\|\varphi\|_{C^3(\overline{\Omega}_+)},\|A^{ij}\|_{C^2(\overline{\Omega}_+)}$,
$\|B^l\|_{C^1(\overline{\Omega}_+)}$\quad $(i,j,l=1,2)$ and
$\|C\|_{C^1(\overline{\Omega}_+)}$.
\end{theorem}
\begin{proof}
By Proposition \ref{skeil} and Lemma \ref{pinl}, we immediately
obtain
\begin{equation}\label{contract}
\|u^\varepsilon\|_{H^1(\Omega_+)}
\leq C\|F^\varepsilon\|_{L^2(\Omega_+)}.
\end{equation}
Inequalities \eqref{bluee} and \eqref{contract} yield
\eqref{fikst}.
\end{proof}
\begin{lemma}
There exist a subsequence
$\{u^{\varepsilon'}\}\subset\{u^\varepsilon\}$ of solution of the
BVP \eqref{market}, such that

$(i)$\,\,\,$u^{\varepsilon'}\rightharpoonup u$\,\,\,weakly
in\,\,\,$H^1(\Omega_+)$, as $\varepsilon'\rightarrow0$;

$(ii)$\,\,\,$u^{\varepsilon'}\rightarrow u$\,\,\,strongly
in\,\,\,$L^2(\Omega_+)$, as $\varepsilon'\rightarrow0$;

Furthermore,

$(iii)$\,\,\,we have the estimates
\begin{equation}\label{symbol}
\|u\|_{H^1(\Omega_+)}\leq C\|F\|_{L^2(\Omega_+)},
\end{equation}
where $C$ is a constant;

$(iv)$\,\,\,$u=0$ on $\Gamma$ in the sense of traces, denoted by
$\gamma(u)=0$.
\end{lemma}
\begin{proof}
 Together
\eqref{fikst} with the fact $H^1(\Omega_+)\hookrightarrow
L^2(\Omega_+)$ compactly, it follows the conclusions of $(i)$,
$(ii)$ and $(iii)$.

In addition, from $u^\varepsilon=0$ on $\Gamma$, the conclusion of
$(i)$, and the fact $\gamma:H^1(\Omega_+)\hookrightarrow
L^2(\Gamma)$ is compact map, it follows easily the conclusions of
$(iv)$.
\end{proof}

\subsection{Some properties of solutions of the BVP \eqref{scienze}}
\begin{lemma}\label{shigerenwu}
$\|(\varphi+\varepsilon)u^\varepsilon\|_{H^2(\Omega_+)}
\leq C\|F\|_{L^2(\Omega_+)}.$
\end{lemma}
\begin{proof}
Consider the following BVP:
$$\left\{
\begin{array}{ll}
A^{ij}[(\varphi+\varepsilon)u^\varepsilon]_{\xi_i\xi_j}
=F^\varepsilon+(2A^{ij}\varphi_{\xi_i}-B^j)u_{\xi_j}^\varepsilon
-cu^\varepsilon\quad\hbox{in}\,\,\,\Omega_+,\\[6pt]
(\varphi+\varepsilon)u^\varepsilon=0\quad\hbox{on}\,\,\,\Gamma.
\end{array}
\right.$$ In the sequel, observe that $\|F^\varepsilon
+(2A^{ij}\varphi_{\xi_i}-B^j)u_{\xi_j}^\varepsilon
-cu^\varepsilon\|_{L^2(\Omega_+)}\leq C\|F\|_{L^2(\Omega_+)}$,
Thus, by employing the regularity theory of second order elliptic
type equation, we obtain the claim.
\end{proof}
Without loss of generality, we may assume $d_0=1$. So
$D_{d_0}^+=D_1^+=[-\pi,\pi]\times[0,1]$. $\eta\in
\mathrm{C}^{\infty}([0,1])$, is defined by
$$\eta=\eta(y)=\left\{\begin{array}{c}
1\quad\hbox{for}\,\,\,0\leq y\leq\frac{1}{2}, \\
0\quad\hbox{for}\,\,\,\frac{2}{3}\leq y\leq 1,
\end{array}\right.
\quad\quad\hbox{and}\quad 0\leq\eta\leq 1.$$
\begin{lemma}\label{fern}
The solution $u^\varepsilon$ of $\pounds^\varepsilon
u=f^\varepsilon$ with $u(x,0)=0$, satisfies that
$$\lim_{\varepsilon\rightarrow 0}
\|\varepsilon u^\varepsilon_y(\cdot,0)\|_{L^2([-\pi,\pi])}=0.$$
\end{lemma}
\begin{proof}
Note that $\pounds^\varepsilon u^\varepsilon=f^\varepsilon$
implies
$$(y+\varepsilon)u^\varepsilon_{yy}
=f^\varepsilon -\omega(y+\varepsilon)u^\varepsilon_{xx}
-au^\varepsilon_x -bu^\varepsilon_y +cu^\varepsilon.$$ Therefore, we
obtain easily
\[(y+\varepsilon)(\eta u^\varepsilon_y)_y
=(y+\varepsilon)\eta_yu^\varepsilon_y
+\eta(y+\varepsilon)u^\varepsilon_{yy}=h^\varepsilon,\] where
$h^\varepsilon:=\eta[f^\varepsilon
-\omega\eta(y+\varepsilon)u^\varepsilon_{xx} -a\eta
u^\varepsilon_x-c\eta u^\varepsilon] +[(y+\varepsilon)\eta_y-\eta
b]u^\varepsilon_y$. So
$$(\eta u^\varepsilon_y)_y=\frac{1}{y+\varepsilon}h^\varepsilon.$$
It follows that
$$\eta u^\varepsilon_y
=-\int_y^1\frac{1}{\tau+\varepsilon}h^\varepsilon(x,\tau)d\tau$$
and
$$\eta(y+\varepsilon)u^\varepsilon_y
=-(y+\varepsilon)\int_y^1\frac{1}{\tau+\varepsilon}h^\varepsilon(x,\tau)d\tau$$
From this inequality we get
$$\aligned
\|\eta(y+\varepsilon)u^\varepsilon_y\|^2_{L^2([-\pi,\pi])}
&=\int_\pi^\pi(y+\varepsilon)^2
\big(\int_y^1\frac{1}{\tau+\varepsilon}h^\varepsilon(x,\tau)d\tau\big)^2 dx\\
&\leq\int_\pi^\pi(y+\varepsilon)^2
\big(\int_y^1\frac{1}{(\tau+\varepsilon)^2}d\tau\big)
\big(\int_y^1|h^\varepsilon(x,\tau)|^2d\tau\big) dx\\
&=(y+\varepsilon)^2
\big(\frac{1}{y+\varepsilon}-\frac{1}{1+\varepsilon}\big)
\int_\pi^\pi\int_y^1|h^\varepsilon(x,\tau)|^2d\tau dx\\
&\leq(y+\varepsilon)
\int_\pi^\pi\int_y^1|h^\varepsilon(x,\tau)|^2d\tau dx\\
&\leq(y+\varepsilon)\|h^\varepsilon\|_{L^2(D_1^+)}^2
\endaligned$$
On the other hand, by a direct calculation and together
\eqref{fikst} with Lemma \ref{shigerenwu}, we get
$$\|h^\varepsilon\|_{L^2(D_1^+)} \leq C\|F\|_{L^2(\Omega_+)}.$$
From the previous two results we deduce
$$\|\eta(y+\varepsilon)u^\varepsilon_y\|_{L^2([-\pi,\pi])}
\leq C(y+\varepsilon)^\frac{1}{2}\|F\|_{L^2(\Omega_+)}.$$ In
particular, for $y=0$, with $\eta(0)=1$, the above inequality yields
$$\|\varepsilon u^\varepsilon_y(\cdot,0)\|_{L^2([-\pi,\pi])}
\leq C\varepsilon^\frac{1}{2}\|F\|_{L^2(\Omega_+)}
\rightarrow 0,\,\, \hbox{as}\,\,\,\varepsilon\rightarrow 0.$$
\end{proof}
\begin{corollary}\label{boundaryproperty}
The solution $u^\varepsilon$ of the BVP \eqref{scienze} satisfies
that
$$\lim_{\varepsilon\rightarrow 0}
\|\varepsilon u^\varepsilon_{\xi_i}\|_{L^2(\Gamma)}=0.$$
\end{corollary}
\begin{proof}
The solution $u^\varepsilon$ of the equation $\pounds^\varepsilon
u=f^\varepsilon$ with $u(x,0)=0$, satisfies that
$u^\varepsilon_x=0$. So it is apparent that $\|\varepsilon
u^\varepsilon_x(\cdot,0)\|_{L^2([-\pi,\pi])}=0.$ By Lemma \ref
{fern} and the transformation of variable, we easily prove the
conclusion of this Corollary.
\end{proof}
\begin{remark}\label{kuhnel}
Under the same condition, we may estimate the solution of the BVP
\eqref{radkeic} in the case $\Omega_-$. By the similar procedure and
argument, it will lead to the completely same results as the case
$\Omega_+$. But we do not state these corresponding results, as well
as omit the proof of these results at here.
\end{remark}

\section{The well-posedness of the BVP \eqref{plenum}.}
In this section, we will discuss the well-posedness of BVP
\eqref{mesi}, i.e. to prove Theorem \ref{lanoitanretni}. For this
purpose, we shall also need the following lemma:
\begin{lemma}[see \cite{Sh}]\label{pulse}
Suppose that $u_+\in H^1(\Omega_+)\,,\,\,u_-\in H^1(\Omega_-)$\,,
and their traces on $\Gamma$ is equal, i.e. $\gamma u_+=\gamma
u_-$\,. A function $u\in L^2(\Omega)$ is defined by
$$u=\left\{
\begin{array}{ll}
u_+\quad\hbox{in}\,\,\,\Omega_+, \\
u_-\quad\hbox{in}\,\,\,\Omega_-.
\end{array}
\right.$$
Then
$u\in H^1(\Omega).$
\end{lemma}
\begin{theorem}\label{wudinghe}
Under the conditions of Theorem \ref{lanoitanretni}. Then there
exists a unique $H^1$ weak solution $u$ of the BVP \eqref{plenum} in
domain $\overline{\Omega}$. And $u$ satisfies
\begin{equation}\label{lingzhou}
\|u\|_{H^1(\Omega)}\leq C\,\|F\|_{L^2(\Omega)},
\end{equation}
where $C$ is a constant depending only on $\Gamma$,
$\|\varphi\|_{C^3(\overline{\Omega})}$,
$\|A^{ij}\|_{C^2(\overline{\Omega})}$,
$\|B^l\|_{C^1(\overline{\Omega})}$ $(i,j,l=1,2)$ and
$\|C\|_{C^1(\overline{\Omega})}$.
\end{theorem}
\begin{proof}
Recall that $u^{\varepsilon'}$ and $u^{(-\varepsilon')}$ are the
solutions of the BVP \eqref{scienze} and  \eqref{oleinik}
respectively. For the convenience, we use $\overline{u}$ and
$\underline{u}$ to express respectively the limits of
$u^{\varepsilon'}$ and $u^{(-\varepsilon')}$, i.e.

(i)\,\,\,$u^{\varepsilon'}\rightharpoonup \overline{u}$\,\,\,
weakly in\,\,\,$H^1(\Omega_+)$,
as $\varepsilon'\rightarrow0$;

(ii)\,\,\,$u^{(-\varepsilon')}\rightharpoonup \underline{u}$\,\,\,
weakly in\,\,\,$H^1(\Omega_+)$,
as $\varepsilon'\rightarrow0$.

A function $u\in L^2(\Omega)$ is defined by
$$u=\left\{
\begin{array}{ll}
  \overline{u}\quad\hbox{in}\,\,\,\Omega_+, \\
  \underline{u}\quad\hbox{in}\,\,\,\Omega_-.
\end{array}
\right.$$
By Theorem \ref{ritein}, Remark \ref{kuhnel}, Lemma \ref{pulse},
and $\gamma(\overline{u})=\gamma(\underline{u})=\underline{u}\big|_{\partial\Omega}=0$,
we obtain that $u\in H^1(\Omega)$, and $u=0$\quad on\,\,$\Gamma\cup\Omega$,
in the sense of traces.

On the other hand, from $L^\varepsilon u^\varepsilon=F^\varepsilon$
and $L^{(-\varepsilon)} u^{(-\varepsilon)}=F^\varepsilon$,
we get
$$(L^\varepsilon u^\varepsilon,v)=(F^\varepsilon,v),
\quad \hbox{for all}\quad v\in C_0^1(\Omega),$$
and
$$(L^{(-\varepsilon)} u^{(-\varepsilon)},v)=(F^\varepsilon,v),
\quad \hbox{for all}\quad v\in C_0^1(\Omega).$$ Thus, by Green
formula, $u^\varepsilon=0$\,\,\,on\,\,\,$\Gamma$, and
$u^{(-\varepsilon)}=0$\,\,\,on\,\,\,$\Gamma\cup\partial\Omega$, we
obtain that the following identities
$$\int_{\Omega_+} [-u^\varepsilon_{\xi_i}(A^{ij}\varphi v)_{\xi_j}
+(\varphi Cu^\varepsilon+B^lu^\varepsilon_{\xi_l})v]d\xi
-\int_\Gamma \varepsilon u^\varepsilon_{\xi_i}A^{ij}vn_jds
=\int_{\Omega_+}Fvd\xi$$ and
$$\int_{\Omega_-} [-u^{(-\varepsilon)}_{\xi_i}(A^{ij}\varphi v)_{\xi_j}
+(\varphi Cu^{(-\varepsilon)}+B^lu^{(-\varepsilon)}_{\xi_l})v]d\xi
-\int_\Gamma \varepsilon u^{(-\varepsilon)}_{\xi_i}A^{ij}vn_jds
=\int_{\Omega_-}Fvd\xi$$ hold for all $v\in C_0^1(\Omega)$, where
$$\vec{n}=
(n_1,n_2)= (\frac{\varphi_{\xi_1}}{|\nabla\varphi|},
\frac{\varphi_{\xi_2}}{|\nabla\varphi|})\big|_{\Gamma}\neq 0$$ is
the unit inward normal direction of $\partial\Omega_+$, Meanwhile
it is also the unit outward normal direction of
$\partial\Omega_-$.

Letting $\varepsilon$ tends to 0,
by the above (i), (ii), Corollary \ref{boundaryproperty} and Remark \ref{kuhnel},
we get
$$\int_{\Omega_+} [-\overline{u}_{\xi_i}(A^{ij}\varphi v)_{\xi_j}
+\varphi C\overline{u} v
+B^l\overline{u}_{\xi_l}v]d\xi
=\int_{\Omega_+}Fvd\xi, \quad \hbox{for all}\,\,v\in C_0^1(\Omega);$$
and
$$\int_{\Omega_-} [-\underline{u}_{\xi_i}(A^{ij}\varphi v)_{\xi_j}
+\varphi C\underline{u} v
+B^l\underline{u}_{\xi_l}v]d\xi
=\int_{\Omega_-}Fvd\xi, \quad \hbox{for all}\,\,v\in C_0^1(\Omega).$$

Consequently, $u$ satisfies that
$$\int_\Omega [-u_{\xi_i}(A^{ij}\varphi v)_{\xi_j}
+\varphi Cu v +B^lu_{\xi_l}v]d\xi =\int_\Omega Fvd\xi \quad
\hbox{for all}\,\,v\in C_0^1(\Omega).$$ Since $C_0^1(\Omega)$ is
dense in $H_0^1(\Omega)$. So, we have thus the existence of $H^1$
weak solution of the BVP \eqref{plenum} on $\Omega$ in a very simple
manner. In addition, \eqref{lingzhou} will be followed by applying
\eqref{symbol} and Remark \ref{kuhnel}.

Finally, together Lemma \ref{weiyixing} with Remark \ref{kuhnel}, we
obtain easily the uniqueness of $H^1$ solution to the BVP
\eqref{plenum} in $\Omega$.
\end{proof}
According to subsection 2.1, we see easily that Theorem
\ref{lanoitanretni} is a trivial Corollary of Theorem
\ref{wudinghe}. So far, we finish the proof for the main result of
this paper.

\subsection*{Acknowledgment}
The results of this paper were obtained during author's Ph.D.
studies at Fudan University and are also contained in his
thesis~\cite{YueHe} with the similar title. The author wishes take
this opportunity to express deep gratitude to his supervisor
Professor JiaXing Hong whose continuous guidance, support, interest
and encourage were crucial for the successful completion of the
present work. The author also thanks Professor Qing Han for his nice
suggestion and enthusiastic help.

% ------------------------------------------------------------------------

\end{document}